\documentclass[11pt]{article}
\usepackage{mathrsfs}
\usepackage{amssymb}
\usepackage{amsmath}
\usepackage[all]{xy}
\newtheorem{Lemma}{Lemma}[section]
\newtheorem{Proposition}{Proposition}[section]
\newtheorem{Theorem}{Theorem}[section]
\newtheorem{Definition}{Definition}[section]

\def\Im{\mathop{\rm Im}\nolimits}
\def\Ker{\mathop{\rm Ker}\nolimits}
\def\Coker{\mathop{\rm Coker}\nolimits}
\def\Tr{\mathop{\rm Tr}\nolimits}

\def\Ext{\mathop{\rm Ext}\nolimits}
\def\Hom{\mathop{\rm Hom}\nolimits}

\title{\Large \bf When Are IG-projective Modules Projective?
}
\author{Rong Luo \thanks{\it E-mail:luorong@swjtu.edu.cn}
\\
{\small \it Collogue of Mathematics, Southwest Jiaotong University,
Chengdu 610031, P. R. China}\\
Dongmei Jian\\
{\small \it College of Math and Software Science, Sichuan Normal University, Chengdu 610066,P. R. China}}
\usepackage{amssymb}
\date{}
\begin{document}
\baselineskip=18pt \maketitle

\begin{abstract} This paper concerns when a finitely generated IG-projective module is projective over commutative Noetherian local rings.
We prove that a
finitely generated IG-projective module is projective if and only if it is selforthogonal.
\end{abstract}

\vspace{0.5cm}

\section{Introduction}
Unless stated otherwise, all rings discussed in this paper are commutative Noetherian local rings, and all modules
are finitely generated. Let $R$ be a commutative Noetherian ring. We use mod$R$ to denote the category of finitely generated $R$-modules.
As a common generalization of the notion of projective modules, Auslander and Bridger
in [AuB] introduced the notion of finitely generated modules of Gorenstein dimension $0$.
Such modules are called Gorenstein projective, following Enochs and Jenda's terminology
in [EJ], which are defined as follows:

\begin{Definition}
{\it An $R$-module $M$ is said to be Gorenstein
projective ($G$- projective, for short) if there exists
an exact sequence of projective modules
$$\xymatrix{{\bf
P}=\cdots\ar[r]&P_1\ar[r]&P_0\ar[r]&P_{-1}\ar[r]&P_{-2}\ar[r]&\cdots }$$
such that $\Hom_R({\bf P}, R)$ is exact and  $M\cong$ $\Im$($P_0\to P_{-1}$).

The exact sequence ${\bf P}$ is called a complete projective
resolution of $M$.}
\end{Definition}
We denote $G(R)$ as the full subcategory of mod$R$ consisting of all Gorenstein projective modules. It is well known that a projective module is Gorenstein projective. It is natural to ask when are the Gorenstein projective modules projective. Our guess is that the
Gorenstein projective module is projective if and only if it is self-orthogonal.  In [LH], it is proved that this conjecture is true if
$R$ is a ring with radical square zero.

\begin{Definition}
An indecomposable $R$-module $M$ is said to be \emph{IG-projective} if  it is G-projective and admits either an
irreducible epimorphism $P\to M$ or an irreducible monomorphism $M\to P$, with
$P$ being a projective module. A (possibly decomposable) module is IG-projective if it is a direct sum of indecomposable IG-projectives.
\end{Definition}

This notion was introduced by Luo[L], who also prove that if, over such an Artin local algebra $R$ with the simple IG-projective module,
then $1$-self-orthogonal modules are projective.

In this paper, one sees  the isomorphisms as irreducible morphisms. Thus,  the projective modules are IG-projective.
The main purpose of this paper is to prove that this conjecture is also true for IG-projective modules if $R$ is
a commutative Noetherian local ring, which is the following theorem.

\begin{Theorem}\label{Thii}
For a commutative Noetherian local ring, a
finitely generated IG-projective module is projective if and only if it is selforthogonal.
\end{Theorem}

In the next section, we start by recalling the definitions of Gorenstein dimension  and  approximation of a module,
give several preliminary lemmas involving their properties.

\vspace{0.5cm}

\section{Preliminaries}
In this section, we provide some background material. Throughout this
section, let $(R,m,k)$ be a commutative Noetherian local ring with the maximal ideal $m$ and the field $k$.
 The starting point is a definition of $G$-dimension, introduced by
Holm[H]

\begin{Definition}
Let $M$ be an $R$-module. If n is a non-negative integer such that there
is an exact sequence$$0\to G_n\to G_{n-1}\to \cdots\to G_1\to G_0\to M\to 0$$
of $R$-modules with $G_i\in G(R)$ for every $i=0,1,\cdots,n,$ then we say that
$M$ has $G$-dimension at most $n$, and write $G-\dim_RM\leq n$. If such an
integer $n$ does not exist, then we say that M has infinite $G$-dimension,
and write $G-\dim_RM=\infty$.
\end{Definition}

Recall $R$ is called Gorenstein if self-injective dimension of $R$ is finite. The next three lemmas are the properties of
$G$-dimension, the proofs are seen in [Ch] and [Ta]
\begin{Lemma}
 Let $0\to L\to M\to N\to 0$ be a short exact sequence of R-modules. If two of $L,M,N$
have finite $G$-dimension, then so does the third.

\end{Lemma}

\begin{Lemma} The following conditions are equivalent:

$(1)$ $R$ is Gorenstein;

$(2)$ $G-\dim_RM <\infty$ for any $R$-module $M$;

$(3)$ $G-\dim_Rk<\infty$.
\end{Lemma}

\begin{Lemma}\label{b1}
Suppose that there is a direct sum decomposition $m=I\oplus J$ where $I,J$
are non-zero ideals and $G$-dim$_RI$ is finite. Then $R$ is a Gorenstein local ring of dimension
one.
\end{Lemma}

Next, the notion of a approximation of a module is introduced by Auslander and Reiten [AuR].

\begin{Definition} Let $\mathscr{X}$ be a
full subcategory of {\rm mod}$R$ and $\phi: X\to M$ be a
homomorphism from $X\in\mathscr{X}$ to $M\in$ mod$R$. We call
$\phi$ a right $\mathscr{X}$-approximation of $M$ if for any
homomorphism $\phi^\prime:X^\prime\to M$ with
$X^\prime\in\mathscr{X}$ there exists a homomorphism $f:X^\prime\to
X$ such that $\phi^\prime=\phi f$.
\end{Definition}

Let $ P_1\to P_0\to M\to 0$
be a presentation with $P_i$ projective $R$-modules. We write $f^*$ for $\Hom_R(f,R)$, $(-)^*$ for $\Hom_R(-,R)$
 and recall that the $R$-module  $\Coker f^*$
is called \emph{the transpose} of $M$, and denote as $\Tr M$; this is well-defined up to projective summands. Here we state
an exact sequence  and isomorphism of functors for later use. For the proofs, we refer to [VM] and [AF].

\begin{Lemma} For any $M\in$
{\rm mod}$R$, there exists an exact sequence of functors from
{\rm mod}$R$ to itself: $$0\to\Ext_R^1(\Tr M,-)\to
M\otimes_R-\stackrel{\lambda(-)}{\longrightarrow}\Hom_R(M^*,-)\to\Ext_R^2(\Tr
M,-)\to 0.$$
\end{Lemma}

\begin{Lemma}

For any $M\in$
{\rm mod}$R$, there exist isomorphisms of functors from
{\rm mod}$R$ to itself:
$$(M\otimes -)^*\cong \Hom_R(M,(-)^*)\cong\Hom_R(-,M^*)$$

\end{Lemma}

\vspace{0.5mm}
\section{The main results}

In this section, let $(R,m,k)$ be a commutative Noetherian local ring with the maximal ideal $m$ and the field
$k$, we begin with introducing a proposition, which plays a crucial role in this
section. Put $D(-)=\Hom_R(-,E(R/J))$ where $J$ is the  Jacobson radical of $R$
and $E(R/J)$ is the injective envelope of $R/J$.

\begin{Proposition}\label{a1}
If $(R,m,k)$ is a local ring such that G-dim$Dm$ is finite,, then there exists an exact sequence $$0\to L\to X\to Dk\to 0$$ with $X$ in $G(R)$ such that

$(1)$ the morphism $X\to Dk$ is a $G(R)$-approximation of $Dk$ and $\Ext_R^i(G(R),L)=0$ for any $i\geq 1$;

$(2)$ the sequence $0\to \Hom_R(Dk,G)\to \Hom_R(X,G)\to \Hom_R(L,G)\to 0$ is exact for any Gorenstein projective $R$-module $G$ that is not projective.
\end{Proposition}


\noindent{\it Proof.} Applying the functor $D(-)$ to the exact sequence $0\to m\to R \to k\to 0$, we
have $0\to Dk\to DR\to Dm \to 0$. Let the morphism $g:Q\to DR$ be a projective cover of $DR$ with the projective module $Q$.
Consider a pull-back diagram of the morphisms $Dk\to DR$ and $Q\to DR$:
$$\xymatrix{
&&0\ar[d]&0\ar[d]\\
0\ar[r]&\Ker g\ar[r]\ar@{=}[d]&Y\ar[r]\ar[d]&Dk\ar[r]\ar[d]&0\\
0\ar[r]&\Ker g\ar[r]&Q\ar[r]\ar[d]&DR\ar[r]\ar[d]&0\\
&&Dm\ar@{=}[r]\ar[d]&Dm\ar[d]\\
&&0&0}$$
then the sequence $0\to \Ker g\to Q\to DR\to 0$ is exact. This induces that $\Ext_R^i(G(R),\Ker g)=0$ for $i>0$.
Since G-dim$_RDm$ is finite, by lemma $2.1$, so is $Y$. We consider the strict $G(R)$-resolution of $Y$, say $0\to P_s\to P_{s-1}
\to \cdots\to P_1\to X\to Y\to 0$ with all the $P_i$ being projective and $X$ belonging to $G(R)$.
Consider the pullback of the morphisms $X\to Y$ and $\Ker g\to Y$
$$\xymatrix{
&&&&&0\ar[d]&0\ar[d]\\
0\ar[r]&P_s\ar@{=}[d]\ar[r]&P_{s-1}\ar@{=}[d]\ar[r]&\cdots\ar[r]&P_1\ar@{=}[d]\ar[r]&L\ar[r]\ar[d]&\Ker g\ar[r]\ar[d]&0\\
0\ar[r]&P_s\ar[r]&P_{s-1}\ar[r]&\cdots\ar[r]&P_1\ar[r]&X\ar[r]\ar[d]&Y\ar[r]\ar[d]&0\\
&&&&&Dk\ar[d]\ar@{=}[r]&Dk\ar[d]\\
&&&&&0&0}$$
then the long sequence
$$0\to P_s\to P_{s-1}\to \cdots\to P_1\to L\to \Ker g\to 0$$
is exact. Therefore $\Ext_R^i(G(R),\Ker g)=0$ tells us that $\Ext_R^i(G(R), L)=0$  for $i>0$.
Thus  $X\in G(R)$ implies that the exact $0\to L\to X\to Dk\to $ is a
$G(R)$-approximation of $Dk$. This completes the proof of $(1)$.

Next to prove $(2)$.  Let $G$ be an any indecomposable $G$-projective $R$-module that is not projective. We take $\lambda_{G^*}(-)$ to be the morphism $\lambda_{G^*}(-): G^*\otimes _R -\to \Hom_R(G^{**}, -)$ by $\lambda_{G^*}(-)(a\times -)(f)=f(a)\cdot -$
for any $a\in {G^*}$, $f\in G^{**}$. Note $\Tr
G^*\in G(R)$ and the $G(R)$-approximation of $Dk$ $g:X\to Dk$, we have
$$\Ker \lambda_{G^*}(L)=\Ext_R^1(\Tr_R G^*, L)=0
{\rm~~and~~} \Coker \lambda_{G^*}(L)=\Ext_R^2(\Tr_R
G^*,L)=0$$  By the lemma $2.4$, this means that $\lambda_{G^*}(L)$ is
an isomorphism. Hence the composite map
$\lambda_{G^*}(X)\cdot({G^*}\otimes _R\theta)={\rm
Hom}_R(G^{**},\theta)\cdot \lambda_{G^*}(L)$ is injective, and
so is the map ${G^*}\otimes _R\theta$. Thus we have the
following commutative diagram
$$\xymatrix{
0\ar[r]&{G^*}\otimes_R
L\ar[rr]^{{G^*}\otimes_R\theta}\ar[d]^{\lambda
_{G^*}(L)}_\cong&&{G^*}\otimes_R X
\ar[rr]^{{G^*}\otimes_R \pi}\ar[d]^{\lambda _{G^*}(X)}&&{G^*}\otimes_R Dk\ar[r]\ar[d]^{\lambda_{G^*}(Dk)}&0\\
0\ar[r]&{\rm Hom}_R (G^{**},L)\ar[rr]^{{\rm Hom}_R
(G^*,\theta)}&&{\rm Hom}_R (G^{**},X)\ar[rr]^{{\rm
Hom}_R (G^{**},\pi)}&&{\rm Hom}_R (G^{**},Dk)\ar[r]&0}$$
with exact rows. Since $G\cong G^{**}$ is a non-projective indecomposable
module, we have $G^*\otimes Dk\to \Hom_R(G^{**},Dk)$ is zero. That is,
$G^*\otimes \theta$ is split and  we have the exact sequence $0\to
(G\otimes Dk)^*\to (G\otimes X)^*\to (G\otimes L)^*\to 0$. Note from the
lemma 2.5, we get the following commutative diagram
$$\xymatrix{
0\ar[r]&(G^*\otimes Dk)^*\ar[r]\ar[d]^\cong&(G^*\otimes
X)^*\ar[r]\ar[d]^\cong&(G^*\otimes L)^*\ar[r]\ar[d]^\cong &0\\
0\ar[r]&\Hom_R(G^*, (Dk)^*)\ar[r]\ar[d]^\cong&\Hom_R(G^*,
X^*)\ar[r]\ar[d]^\cong&\Hom_R(G^*, L^*)\ar[r]\ar[d]^\cong&0\\
0\ar[r]&\Hom_R(Dk, G^{**})\ar[r]&\Hom_R(X,
G^{**})\ar[r]&\Hom_R(L, G^{**})\ar[r]&0}$$ That is,
$$0\to
\Hom_R(Dk, G)\to \Hom_R(X, G)\to \Hom_R(L, G)\to
0$$for any non-projective module $G$ in $G(R)$ .\hfill{$\square$}

\vspace{1.8mm}

Let $M$ be in $G(R)$. We denote $\Omega ^1(M)$ to be the $1$th syzygy module of $M$. By the definition of
Gorenstein projective module, $\Omega^1(M)$ is in $G(R)$.

\begin{Proposition}\label{a2}  If $(R,m,k)$ is a local ring such that G-dim$Dm$ is finite,
then any indecomposable $IG$-projective $R$-module $M$ satisfying $\Ext_R^i(M,M)=0$ for $i\geq 1$ is projective.

\end{Proposition}

\noindent{\it Proof.}  Assume that
$M$ is non-projective. We want to derive a contradiction. Since $M$ is Irre-Gorenstein
projective, there exists the irreducible morphism $f:P\to M$ or $h:M\to P$ with a projective module $P$.

$(1)$ If such an $f$ exists, then we take a non-split exact sequence $0\to k\to E^\prime\to M\to 0$. Since $f$ is irreducible, it follows that $E^\prime\cong P\oplus E_1$ and the following diagram
$$\xymatrix{
&0\ar[d]&0\ar[d]\\
0\ar[r]&K\ar[r]\ar[d]&P\ar[r]^f\ar[d]&M\ar[r]\ar@{=}[d]&0\\
0\ar[r]&k\ar[r]\ar[d]&P\oplus
E_1\ar[r]\ar[d]&M\ar[r]&0\\
&E_1\ar[d]\ar@{=}[r]&E_1\ar[d]\\
&0&0}$$
is commutative.
If $K=0$, then $M$ is projective. If $E_1=0$, there is the exact sequence $0\to k\to P\to M\to 0$.
Since $\Ext_R^i(M, M)=0$ for $i>0$, we have $\Ext_R^2(M, k)=0$. That is, pd$_RM$ is finite. Hence, $M$ is projective.

$(2)$ Assume that $h$ exists. Since $G$-dim$Dm$ is finite, by proposition \ref{a1}, there exists
a short exact sequence:
 $$0\to L\to X\to Dk\to 0$$ of $R$-modules such that $X\to Dk\to 0$ is a $G(R)$-approximation of $Dk$ and $\Ext_R^i(G(R),L)=0$ for $i>0$.
Take a non-split exact sequence $0\to M\to E\to Dk\to 0$ in
$\Ext_R^1(Dk, M)$, we have the pullback diagram:
$$\xymatrix{ &&0\ar[d]&0\ar[d]\\
&&M\ar[d]\ar@{=}[r]&M\ar[d]
\\
0\ar[r]&L\ar[r]\ar@{=}[d]&Q\ar[r]\ar[d]&E\ar[r]\ar[d]&0\\
0\ar[r]&L\ar[r]&X\ar[r]\ar[d]&Dk\ar[r]\ar[d]&0\\
&&0&0}
$$
with $Q\in G(R)$.

a) If $0\to M\to Q\to X\to 0$ is split, by the (2) of proposition \ref{a1}, there exists the following commutative diagram
$$\xymatrix{&0\ar[d]&0\ar[d]\\
0\ar[r]&\Hom_R(Dk,M)\ar[r]\ar[d]&\Hom_R(X,M)\ar[r]\ar[d]&\Hom_R(L,M)\ar[r]\ar@{=}[d]&0\\
0\ar[r]&\Hom_R(E,M)\ar[r]\ar[d]&\Hom_R(Q,M)\ar[r]\ar[d]&\Hom_R(L,M)\ar[r]&0\\
&\Hom_R(M,M)\ar@{=}[r]&\Hom_R(M,M)\ar[d]\\
&&0}$$
This induces that $0\to \Hom_R(Dk,M)\to \Hom_R(E,M)\to\Hom_R(M,M)\to 0$ is exact. So we have the
exact sequence $0\to M\to E\to Dk\to 0$ is split. This is contradicted with it being non-split.

b) The next, let $0\to M\to Q\to X\to 0$ be non-split. Note from $M\in G(R)$, we take a short exact sequence
$$0\to M\to P\to M_0\to 0$$ with $M_0\in G(R)$. Since the monomorphism $M\to P$ is irreducible and
$X$ is in $G(R)$, there is  the following commutative diagram
$$\xymatrix{0\ar[r]&M\ar[r]\ar@{=}[d]&Q\ar[r]\ar[d]^\theta&X\ar[r]\ar[d]&0\\
0\ar[r]&M\ar[r]&P\ar[r]&M_0\ar[r]&0}$$
where the morphism $\theta$ is split epimorphic. That is, there is an exact sequence $0\to Q_0\to X\to M_0\to 0$ with $\Ker\theta=Q_0$. Since
$L$ is the maximal submodule of $X$, we get the following commutative diagram
$$\xymatrix{&0\ar[d]&0\ar[d]\\
&Q_0\ar[d]\ar@{=}[r]&Q_0\ar[d]\\
0\ar[r]&L\ar[r]\ar[d]&X\ar[r]\ar[d]&Dk\ar[r]\ar@{=}[d]&0\\
0\ar[r]&N\ar[r]\ar[d]&M_0\ar[r]\ar[d]&Dk\ar[r]&0\\
&0&0}$$
 It follows that, by the $(1)$ of proposition \ref{a1}, an exact sequence
$0\to \Hom_R(\Omega^1(M),L)\to \Hom_R(\Omega^1(M),X)\to \Hom_R(\Omega^1(M),Dk)\to 0$. This implies the sequence
$$0\to \Hom_R(\Omega^1(M),N)\to \Hom_R(\Omega^1(M),M_0)\to \Hom_R(\Omega^1(M),Dk)\to 0\eqno{(*)}$$
is exact.

Consider the  push-out diagram of the morphisms $M_0\to Dk$ and $E\to Dk$:
$$\xymatrix{&&0\ar[d]&0\ar[d]\\
&&M\ar[d]\ar@{=}[r]&M\ar[d]\\
0\ar[r]&N\ar[r]\ar@{=}[d]&Q'\ar[r]\ar[d]&E\ar[r]\ar[d]&0\\
0\ar[r]&N\ar[r]&M_0\ar[r]\ar[d]&Dk\ar[r]\ar[d]&0\\
&&0&0}$$
we have the sequence $\aleph : 0\to M\to Q^\prime\to M_0\to 0$ is exact. If $\aleph$ is split, then this induces a contradiction by repeating  the proceedings of $a)$.  Let $\aleph$ be a non-split exact sequence. Since $M_0$ is in $G(R)$, there exists the following commutative diagram
$$\xymatrix{
0\ar[r]&M\ar[r]\ar@{=}[d]&Q^\prime \ar[r]\ar[d] &M_0\ar[r]\ar[d]&0\\
0\ar[r]&M\ar[r]&P\ar[r]&M_0\ar[r]&0}$$
Since $M\to P$ is irreducible, one have the morphism $Q^\prime\to P$ is split epimorphic. We easily
see that $P\cong Q^\prime$. Thus we obtain a commutative diagram
$$\xymatrix{&&0\ar[d]&0\ar[d]\\
&&M\ar[d]\ar@{=}[r]&M\ar[d]\\
0\ar[r]&N\ar[r]\ar@{=}[d]&P\ar[r]\ar[d]&E\ar[r]\ar[d]&0\\
0\ar[r]&N\ar[r]&M_0\ar[r]\ar[d]&Dk\ar[r]\ar[d]&0\\
&&0&0}$$

\textbf{We claim that  $\Ext_R^2(M, N)=0$}. Since $M$ is  selforthogonal, by the exact sequence $0\to M\to P\to M_0\to 0$, then $\Ext_R^1(M,M_0)=0$.  Note from our claim, by the exact sequence $0\to N\to M_0\to Dk\to 0$, we have $\Ext_R^1(M,Dk)=0$. This is, $M$ is projective.

\textbf{Next to prove our claim. } Since $\Omega^1(M)$ is the $1$th syzygy of $M$, the selforthogonal module $M$ implies that $\Ext_R^1(\Omega^1(M), M)=0$ . Applying the functor $\Hom_R(\Omega^1(M),-)$ to the above diagram,  we get a commutative diagram
$$\xymatrix{&&0\ar[d]&0\ar[d]\\
&&\Hom_R(\Omega^1(M), M)\ar[d]\ar@{=}[r]&\Hom_R(\Omega^1(M),M)\ar[d]\\
0\ar[r]&\Hom_R(\Omega^1(M),N)\ar[r]\ar@{=}[d]&\Hom_R(\Omega^1(M),P)\ar[r]\ar[d]&\Hom_R(\Omega^1(M),E)\ar[d]\\
0\ar[r]&\Hom_R(\Omega^1(M),N)\ar[r]&\Hom_R(\Omega^1(M),M_0)\ar[r]^\delta\ar[d]&\Hom_R(\Omega^1(M),Dk)\\
&&0}$$
   Note from  the exact sequence $(*)$ that $\delta$ is epimorphic. Thus we get an exact sequence
   $$0\to \Hom_R(\Omega^1(M),N)\to \Hom_R(\Omega^1(M),P)\to \Hom_R(\Omega^1(M),E)\to 0$$
 This induces  $\Ext_R^1(\Omega^1(M),N)=0$. The $1$th syzygy $\Omega^1(M)$ tells us that $\Ext_R^2(M,N)=0$.

The results of $(1)$ and $(2)$ contrary to the assumption of the proposition. This contradiction completes the
proof of the proposition.
\hfill{$\square$}

\vspace{0.2cm}
Now, let us prove our main theorem.
\begin{Theorem}
Let $(R, m, k)$ be a commutative Noetherian local ring. An  $IG$-projective  $R$-module $M$ is projective if and only if $M$ is selforthogonal.
\end{Theorem}

\noindent{\it Proof.} Without loss of generality, let $M$ be an indecomposable module. If $R$ be a Gorenstein ring, then G-dim$Dm$ is finite. By the proposition $\ref{a2}$, we have our result.

Let $R$ be a non-Gorenstein ring. Assume that $M$ is an  non-projective module.
We need to derive a contradiction. Since $M$ is Irre-Gorenstein
projective, there exists the irreducible morphism $f:P\to M$ or $h:P\to M$ with a projective module $P$.

$(1)$ If such an $f$ exists, then taking a non-split exact sequence $0\to k\to E\to M\to 0$ and arguing as in the proof $(1)$ of proposition $\ref{a2}$, one deduces that $M$ is projective.

$(2)$ If there exists the irreducible monomorphism $h:M\to P$. Since $M$ is not projective, there exists a maximal submodule $M_1/M$ of
$P/M$. Consider the commutative diagram
$$\xymatrix{&&M_1\ar[rd]^l\\0\ar[r]&M\ar[ru]^{h_1}\ar[rr]^h&&P}$$
the irreducible morphism $h$ implies that $h_1$ is split monomorphic. That is, $M_1=M\oplus H$ and $M_1$ is a maximal submodule of $P$.
Hence, there exists  the following commutative diagram
$$\xymatrix{0\ar[r]&m\ar[r]\ar[d]&R\ar[r]\ar[d]&k\ar[r]\ar@{=}[d]&0\\
0\ar[r]&M\oplus H\ar[r]&P\ar[r]&k\ar[r]&0}$$
By the Schanuel's lemma, we have the isomorphism $m\oplus P\cong R\oplus M\oplus H$. Since $M$ is not projective, $M$ is a summand of $m$.
Since $R$ is a non-Gorenstein ring, it is contradicted with the Lemma $\ref{b1}$.  Hence, $M$ is projective.

The results of $(1)$ and $(2)$ contrary to the assumption of the theorem. This contradiction completes the
proof of the theorem.

\vspace{0.5cm}



\begin{thebibliography}{101}


\bibitem[AF]{A1} F.W. Anderson, K.R. Fuller, Rings and Categories of Modules, second ed., Grad. Texts in Math., vol. 13,Springer-Verlag, Berlin, 1992.


\bibitem[AuB]{A3} M. Auslander and M. Bridger, Stable module theory,
Memoirs Amer. Math. Soc. {\bf 94}, Amer. Math. Soc., Providence,
Rhode Island, 1969.


\bibitem[AuR]{A6} M. Auslander and I. Reiten,  Applications of
contravariantly finite subcategories, Adv. Math. {\bf 86}(1991),
111--152.

\bibitem[Ch]{A6} L.W. Christensen, Gorenstein dimensions. Lecture Notes in Mathematics 1747,
Springer, Berlin, 2000

\bibitem[EJ]{A9} E.E. Enochs and O.M.G. Jenda,  Gorenstein injective
and projective modules, Math. Z. {\bf 220}(1995), 611--633.



\bibitem[H]{A14} H. Holm, Gorenstein homological dimensions, Journal of Pure and Applied Algebra {\bf 189} (2004) 167--193.

\bibitem[LH]{A14} R. Luo and Z. Huang, When are torsionless modules projective?, Journal of Algebra {\bf 320} (2008) 2156--2164.

\bibitem[L]{A14} R. Luo, IG-projective modules, Journal of Pure and Applied Algebra, http://dx.doi.org/10.1016/j.jpaa.2013.05.010.

\bibitem[Ta]{A12} R. Takahashi,  On the category of modules of Gorenstein dimension zero II,
JournalAlgebra {\bf 278} (2004), 402--410.


\bibitem[VM]{A15} Vladimir Masek,  Gorenstein dimesnion and
torsion of modules over commutative Noetherian rings, Comm. Algebra
{\bf 28}(12)(2000), 5783--5811.

\end{thebibliography}
\end{document}